# SEMI-CARTESIAN SQUARES

# AND THE SNAKE LEMMA

JEAN-CLAUDE RAOULT[1]

*À Claude Chevalley, dont the cours of DEA (1965-1966)
m'a fait découvrir les carrés semi-cartesians.*

Résumé: On démontre le lemme du serpent de façon purement catégorique (§ 3). Aucun point n'apparaîtra, ni « points » au sens of Grothendieck ni pseudo-éléments (Guglielmetti & Zaganidis [2009]). En revanche, un fort usage sera fait des carrés semi-cartésiens (§ 2) introduits par Chevalley. Le paragraphe 1 is dévolu à quelques résultats de base sur les catégories abéliennes utilisés par la suite.

Abstract: The snake lemma is proved entirely within category theory (§ 3) without the help of "points with value in..." *à la* Grothendieck nor pseudo-elements (Guglielmetti & Zaganidis [2009]). Instead, we use consistently semi-cartesian squares (§ 2), promoted by Chevalley. Section 1 is devoted to a few basic results on abelian categories, for further use.

This paper is mainly intended to promote the semi-cartesian squares, introduced by Chevalley in a course given at the IHP, and is an example of their flexibility. The first two sections are extracted from this course. The third is a purely categorical proof of the snake lemma.

Categories are supposed to be known: objects and arrows between objects. Arrows are composed associatively and each object X has an identity arrow denoted $1_X$. The *dual* or *opposite* category has same objects but arrows are reversed. The class of arrows from X to Y is denoted by $\mathrm{Hom}_\mathbf{C}(X,Y)$. In a *small category*, arrows form a set (objects also, why?). An *initial* (*final*) object is an objet with a unique arrow to (from) every object. By definition, *monomorphisms* are arrows that are simplifiable from the left ($mu = mv \Rightarrow u = v$) and *epimorphisms* arrows simplifiable from the right ($up = vp \Rightarrow u = v$). Categories will be denoted by bold upper case letters.

A *functor* is a mapping between categories F : $\mathbf{C} \to \mathbf{D}$ compatible with the identities and composition. Given two functors F and G from $\mathbf{C}$ to $\mathbf{D}$, a *natural* (or *functorial*) *morphism* $\varphi : F \to G$ is a family of arrows of $\mathbf{D}$ indexed by the objects of $\mathbf{C}$, such that the squares of figure 1 are commutative for all arrows X → Y of $\mathbf{C}$. If $\mathbf{C}$ is a small category, functors from $\mathbf{C}$ to $\mathbf{D}$ are the objects of a category $F(\mathbf{C}, \mathbf{D})$ the arrows of which are the natural morphisms.

$$\begin{array}{ccc} X & \xrightarrow{f} & Y \\ & & \\ F(X) & \xrightarrow{F(f)} & F(Y) \\ \varphi_X \downarrow & & \downarrow \varphi_Y \\ G(X) & \xrightarrow{G(f)} & G(Y) \end{array}$$

**Figure 1.**

---

1    *20 avenue J. B. Marrou, 63122 Ceyrat,* `jean-claude.raoult@sfr.fr`

## 1. Summary of abelian category

Given a small category **I** of indices, a functor $A : \mathbf{I} \to \mathbf{C}$ can be seen as a commutative diagram of type **I** in **C**. With every object X of **C** is associated a constant diagram $(K_X)_i = X$ and $(K_X)(1_i) = 1_X$. A *projective* (*inductive*) *limit* of a functor $A : \mathbf{I} \to \mathbf{C}$ is a right (left) adjoint of functor $K : X \to K_X$ i.e. $\mathrm{Hom}_{F(I,C)}(K_X, A) \simeq \mathrm{Hom}_C(X, \mathrm{limproj}\, A)$ (resp. $\mathrm{Hom}(A, K_X) \simeq \mathrm{Hom}(\mathrm{limind}\, A, X)$).

*Examples*: 1) $I = \emptyset$: the projective (inductive) limit of the empty set is the final (initial) object.

2) $I = \{1,2\}$: the projective limit of $A_1$, $A_2$ is the product $A_1 \times A_2$ and the inductive limit is the sum $A_1 + A_2$.

3) $I = \{1 \rightrightarrows 2\}$: the projective limit of a double arrow $(u,v) : A_1 \rightrightarrows A_2$ is the *kernel* or *equalizer* of $(u,v)$. The inductive limit is its *cokernel* or *coequalizer*.

4) $I = \{1 \to 0 \leftarrow 2\}$: the projective limit of $A_1 \to A_0 \leftarrow A_2$ is the *cartesian square* or *fiber product* built on these arrows (left square beneath). The inductive limit is got by reversing the arrows: it is a *cocartesian square* or *amalgamated sum* (right square beneath).

$$\begin{array}{ccc} P & \longrightarrow & A_2 \\ \downarrow & & \downarrow \\ A_1 & \longrightarrow & A_0 \end{array} \qquad \begin{array}{ccc} A_0 & \longrightarrow & A_1 \\ \downarrow & & \downarrow \\ A_2 & \longrightarrow & S \end{array}$$

DEFINITION 1. — A *category is abelian if*
*1) it contains a null objet (i. e. initial and final) denoted by* 0;
*2) it accepts finite projective and inductive limits;*
*3) every monomorphism is a kernel and every epimorphism is a cokernel.*

These axioms are due to P. Freyd [1964]. They are preserved by duality.

The null arrow $A \to B$ is the composed arrow $A \to 0 \to B$; the kernel of a single arrow $u : A_1 \to A_2$ is the kernel of $(u,0)$. It is easy to show that every kernel is a monomorphism and that every cokernel is an epimorphism. Condition (3) shows that these two notions coincide and more precisely:

LEMMA 1. — *1) If n is a kernel of an epimorphism q, then q is a cokernel of n.*
*2) If q is a cokernel of a monomorphism n, then n is a kernel of q.*

Case 2 is dual of case 1. For cas 1, $q$ is a cokernel of an arrow $f$: $qf = 0$ and there exists therefore a unique arrow $g$ such that $f = ng$ because $n$ is a kernel of $q$. Suppose an arrow $u$ such that $a = 0$. Then $ung = uf = 0$ and since $q$ is a cokernel of $f$, there exists a unique arrow $v$ such that $u = vq$, QED.

One can deduce the following decomposition of the arrows.

PROPOSITION 1. — *Every arrow can be decomposed into $f = mq$ where m is a monomorphism and q an epimorphism. This decomposition is unique up to a unique isomorphism.*

Soit $m$ the kernel d'une arrow cokernel $p$ of $f$. Then $pf = 0 \Leftrightarrow (\exists! q)\, f = mq$. To show that $q$ is an epimorphism, let us remark:

a) $f$ epimorphism $\Leftrightarrow p = 0 \Leftarrow m$ invertible (because $pm = 0 \Rightarrow pmm^{-1} = p = 0$)

b) Decompose $q$ like $f$: $n$ is a kernel of a cokernel of $q$. Considering a), it is enough to show that $n$ is invertible. Now let $r$ stisfying $rmn = 0$. Then $rmns = rf = 0$. Since $p$ is a cokernel of $f$, $(\exists!t)\, r = tp$. Hence $p$ is a cokernel of $mn$ and from lemma 1, $mn$ is a kernel of $p$. But $m$ is also a kernel of $p$, therefore $n$ is invertible.

c) Let $f = m'q'$ be an other decomposition. Since $p$ is a cokernel of $m'q'$ and $q'$ is an epimorphism, $p$ is a cokernel of $m'$. By lemma 1, $m'$ is a kernel of $p$. Since $m$ is another one, they are isomorphic, QED.

COROLLARY. — *An arrow which is a monomorphism and an epimorphism is an isomorphism.*

The decomposition of $f$ into an epimorphism followed by a monomorphism is unique up to isomorphisms. But $f$ has two such decompositions: $1f = f1$. Therefore it is invertible, QED.

PROPOSITION 2. — *If $q$ is an epimorphism, then $mq$ and $m$ have same cokernel. The converse is true if $m$ is a monomorphism. If $m$ is a monomorphism, then $mq$ and $m$ have same kernel. The converse is true is $m$ is an epimorphism.*

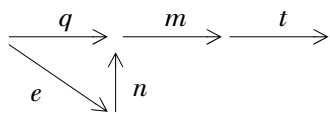

If $q$ is an epimorphism, every arrow $t$ satisfies $tm = 0 \Leftrightarrow tmq = 0$. Hence, $m$ and $mq$ have same cokernel. Conversely, decompose $q$ into an epimorphism $e$ and a monomorphism $n$: The direct part shows that $mne$ and $mn$ have same cokernel. The hypothesis becomes: monomorphisms $m$ and $mn$ have same cokernel. They are therefore kernel of the same arrow; they are isomorphic and $n$ is invertible: $q$ is an epimorphism.

The second assertion is the dual of the first, QED.

PROPOSITION 3. — *If $ba$ is a kernel of $c$ and $b$ is a monomorphism, then $a$ is a kernel of $cb$.*

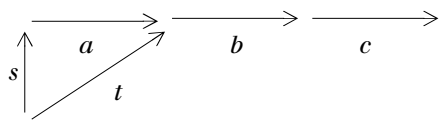

Let $t$ be an arrow with $cbt = 0$. Since $ba$ is a kernel of $c$, there exists a unique arrow $s$ such that $bt = bas$. Since $b$ is a monomorphism, $t = as$, QED.

The following notion is a well-known generalisation of the notion of kernel and cokernel.

DEFINITION 2. — *Let two successive morphisms $f = mq$ and $g = np$, decomposed into epimorphisms followed by monomorphisms. The sequence $(f, g)$ is exact when $m$ is a kernel of $p$.*

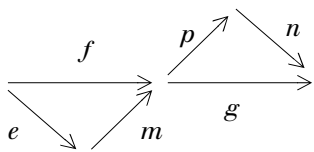

Equivalently (prop. 2), one can require that $m$ be a kernel of $g$, or that $p$ be a cokernel of $m$, or that $p$ be a cokernel of $f$.

Finally, recall that in an abelian category, there exists an isomorphism from the sum to the product and that the insertions $i : A \to A+B$ and $j : B \to A+B$, and the projections $p : A\times B \to A$ and $q : A\times B \to B$ satisfy

$$pi = 1, \qquad qi = 0,$$
$$pj = 0, \qquad qj = 1.$$
$$ip + jq = 1$$

These equalities caracterize the direct sum of A and B, which will be denoted by A + B.

**2. Semi-cartesian squares**

What does one get by composing cartesian and cocartesian squares? Semi-cartesian squares in the following sense.

PROPOSITION 1. — *Let $ca = db$ be a commutative square as in Fig. 2. Let $B + C$ be the direct sum with insertions $(i,j)$ and projections $(p,q)$; construct the fiber product $(P, f, g)$ of $(c,d)$ with the kernel $n$ of*

*cp – dq*, and the amalgamated sum (S,*r*,*s*) *of* (*a*,*c*) *with the cokernel t of ia + jb. Then there exist unique arrows e* : A → P *and m* : S → C *making commutative the obvious triangles of figure 2. Then the following conditions are equivalent:*

*(i) e is an epimorphism,*

*(ii) m is a monomorphism,*

*(iii) the sequence* 0 → P $\xrightarrow{n}$ B+C $\xrightarrow{t}$ S → 0 *is exact,*

*(iv) the sequence* A $\xrightarrow{ia+jb}$ B+C $\xrightarrow{cp-dq}$ D *is exact.*

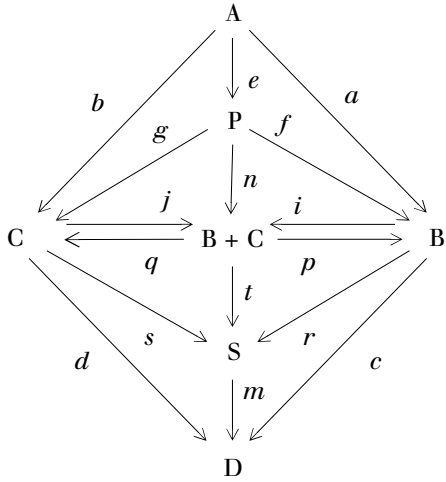

**Figure 2.**

When constructing P and S we defined

$$f = pn, \quad \text{and} \quad r = ti,$$
$$g = qn, \quad\quad\quad\quad s = -tj.$$

So there exists a unique *e* with *a* = *fe* and *b* = *ge*; and a unique *m* with *c* = *mr* and *d* = *ms*. Composing on the right $1_S = ip + jq$ with *ne*, one finds *ne* = *ia* + *jb*. Therefore *t* is a cokernel of *ne*.

(i) ⇔ (iii): For the sequence in (iii) to be exact, it is necessary and sufficient that *t* be a cokernel of *n*, that is, that *e* be an epimorphism.

(iii) ⇔ (ii): Condition (iii) is preserved by duality and is thus is equivalent to (ii) which is dual of (i).

(iii) ⇒ (iv): because *n* is a monomorphism and *q* an epimorphism.

(iv) ⇒ (iii): because the cokernel *t* of *ne* is cokernel of a kernel of *mt* = *cp* – *dq*, that is of *n*, hence (iii), QED.

DEFINITION. — *A commutative square is semi-cartesian if it satisfies the conditions of proposition 1.*

For instance, a cartesian square (*e* is invertible), or a cocartesian square (*m* is invertible), is semi-cartesian. Next is a partial converse in which notations are those of figure 2.

PROPOSITION 2. — *In a semi-cartesian square ca = db, if a is a monomorphism, then d is a monomorphism and the square is cartesian. Si d is an epimorphism, then a is an epimorphism and the square is cocartesian.*

With the notations of Fig. 2, since *a* is a monomorphism, *e* is also a monomorphism. Since it is an epimorphism, it is invertible and the given square is cartesian. Let *k* : N → C be a kernel of *d* and 0 : N → B the null arrow. There exists a unique arrow *h* : N → A such that *k* = *bh* and 0 = *ah*. But *a* is a monomorphism, therefore *h* = 0, hence *k* = 0, QED.

Contrary to arrows, squares will be written in the same order as they are drawn.

PROPOSITION 3. —

*1) Suppose* K *is cocartesian. Then* KL *semi-cartesian* ⇔ L *semi-cartesian.*

*2) Suppose* L *is cartesian. Then* KL *semi-cartesian* ⇔ K *semi-cartesian.*

*3)* K *and* L *semi-cartesians* ⇒ KL *semi-cartesian.*

1) Let (*r*, *s*) be an amalgamated sum of (*c*, *v*) and *m* the unique arrow such that *w* = *mr* and *d* = *ms*. The square *r*(*ca*) = (*sb*)*u* is composed of cocartesian squares and therefore is cocartesian. then KL semi-cartesian ⇔ *m* monomorphism ⇔ L semi-cartesian (see Fig. 3).

SEMI-CARTESIAN SQUARES AND THE SNAKE LEMMA 5

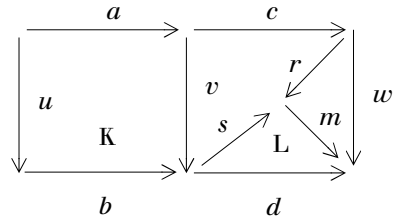

**Figure 3.**

2) Dual of (1).

3) Consider figure 4; set $S = B +_A C$, with a unique monomorphism $n : S \to D$; set $T = S +_B E$. The square ACTE is then cocartesian (composition of cocartesian squares) hence a unique arrow $m : T \to F$ making the diagram commutative. Since L is semi-cartesian and BSTE cocartesian, square SDFT is semi-cartesian from (1) and $m$ is a monomorphism from proposition 2, QED.

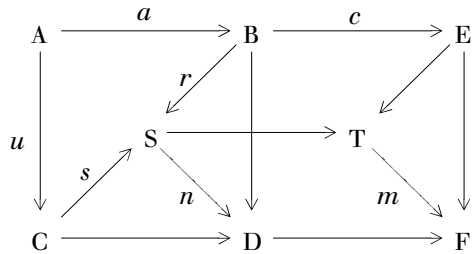

**Figure 4.**

This proposition shows that a semi-cartesian square remains a semi-cartesian square when is removed a cartesian square on the right or a cocartesian square on the left; and also that semi-cartesian squares are got by composing cartesian and cocartesian squares. This is always the case, as shown by the corollary of the following proposition.

PROPOSITION 4. — *Let* KL *be a semi-cartesian square. If* K *is an epimorphism, then* L *is semi-cartesian. If* L *is a monomorphism, then* K *is semi-cartesian.*

Let us prove the first assertion; the second is dual.

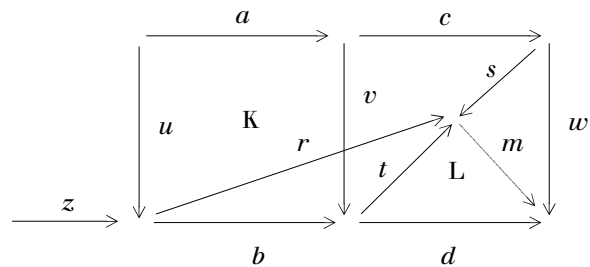

**Figure 5.**

Let $(r,s)$ be an amalgamated sum of $(u,ca)$; since KL is semi-cartesian, there exists a unique monomorphism $m$ with $mr = db$. Since $b$ is an epimorphism, it is a cokernel of some $z$, and $rz = 0$ (compose on the left with $m$ monomorphism). This implies a unique $t$ such that $r = tb$; and $tv = sc$ (compose on the right with epimorphism $a$).

Le square $tv = sc$ is cocartesian: if $xc = yv$, a fortiori $xca = yva = ybu$ and since $(ru = s(ca))$ is a cocartesian square, there exists a unique $n$ such that $x = ns$ and $yb = nr = ntb$. Since $b$ is an epimorphism, one deduces $y = nt$. Since $m$ is a monomorphism, square L is semi-cartesian, QED.

In other terms, in the class of semi-cartesian squares, one can simplify by epimorphisms on the left and by monomorphisms on the right. Beware that the converse is false: a semi-cartesian square (for instance the identity) preceded by an epimorphism is not necessarily semi-cartesian (there exists epimorphisms that are not semi-cartesian).

Corollary. — *Semi-cartesian squares can be decomposed into a cocartesian epimorphism followed by a cartesian monomorphism.*

Decompose the squre into an epimorphism followed by a monomorphism. Proposition 4 ensures that they are also semi-cartesian squares. From proposition 2, the first one is cocartesian and the second one is cartesian.

Proposition 5. — *Consider two successive commutative squares* K *and* L *as in figure 6.*
*1) Suppose that* K *is a kernel of* L. *Then:*
*a) w monomorphism ⇒ K cartesian.*
*b) L semi-cartesian ⇒ u epimorphism.*
*2) Dualy suppose that* L *is a cokernel of* K. *Then:*
*c) u epimorphism ⇒ L cocartesian.*
*d) K semi-cartesian ⇒ w monomorphism.*

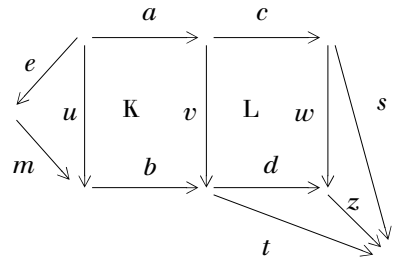

**Figure 6.**

Assertions (a),(b) are dual of (c),(d). Let us show (c) and consider figure 7, in which ($s$, $t$) verifies only $sc = tv$. Then $0 = sca = tva = tbu$ hence $tb = 0$ since $u$ is an epimorphism. Since $d$ is cokernel of $b$, there exists a unique arrow $z$ such that $t = zd$. One checks $s = zw$ by composing with epimorphism $c$ on the right.

Let us show (b). After decomposing L with the help of the above corollary, one may assume that L is a cocartesian epimorphism. Let $u = me$ be the decomposition of $u$ into an epimorphism followed by a monomorphism. If one shows that $d$ is a cokernel of $bm$, then $bm$ will be a kernel of $d$ (cf. lemma 1 § 1) as is $b$, and therefore $m$ will be invertible. Now let $t$ be such that $tbm = 0$. It follows $tbme = 0 = tva$. Since $c$ is a cokernel of $a$, there is a unique arrow $s$ such that $sc = tv$. Since L is cocartesian, there exists a unique arrow $z$ such that $t = zd$ (and $s = zw$), QED.

### 3. The snake lemma

The snake lemma constructs an exact sequence connecting kernels and cokernels.

Proposition 1. — *Suppose two successive squares* K *and* L, *where* L *is semi-cartesian. If* ($a$,$c$) *is exact and $db = 0$, then* ($b$,$d$) *is exact. Dualy, supposing* K *semi-cartesian, then if* ($b$,$d$) *is exact and $ca = 0$, then* ($a$,$c$) *is exact.*

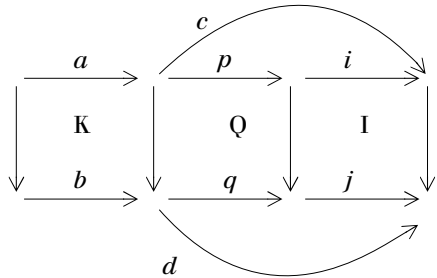

**Figure 7.**

Let Q be the cokernel of K so that L = KI (see Fig. 7). Since $ca = 0$ and $db = 0$, there exist unique arrows $i$ and $j$ such that $c = ip$ and $d = jq$. Then I is semi-cartesian proposition 4 § 2. Since the sequence $(a,c)$ is exact, $i$ is a monomorphism. From proposition 2 § 2, $j$ also is a monomorphism (and I is cartesian): $(b,d)$ is exact, QED.

For each arrow $u$ one selects a kernel arrow of $u$ and denotes its source by $\mathrm{Ker}(u)$. In this way, $\mathrm{Ker}(u)$ becomes a functor.

PROPOSITION 2. — *Kernel functors are left-exact; cokernel functors are right-exact.*

Kernels are (finite) projective limits. Therefore, they commute with projective limits. Dualy, cokernel functors Coker are right-exact.

The following lemma, called the snake lemma, connects these two functors.

LEMMA. — *Given a diagram like Fig. 7, in which i, j, k are kernels of u, v, w, and p, q, r are their cokernels, in which c is a cokernel of a and b is a kernel of d, there exists an arrow δ such that the following sequence is exact:*

$$\mathrm{Ker}(u) \to^s \mathrm{Ker}(v) \to^t \mathrm{Ker}(w) \to^\delta \mathrm{Coker}(u) \to^x \mathrm{Coker}(v) \to^y \mathrm{Coker}(w).$$

Decompose $a = me$ into an epimorphism $e$ followed by a monomorphism $m$ as in figure 7. There exists arrows $i'$ and $u'$ such that $js = mi'$ and $vm = bu'$, because $s$, $m$ and $b$ are the respective kernels of $t$, $c$ and $d$. And since the functor Ker is left exact, $i'$ is a kernel of $u'$. In this way, changing notations, one may assume that $a$ is a kernel of $c$ and dualy that $d$ is a cokernel of $b$.

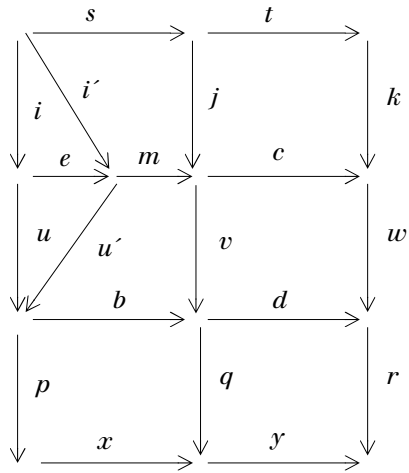

**Figure 8.**

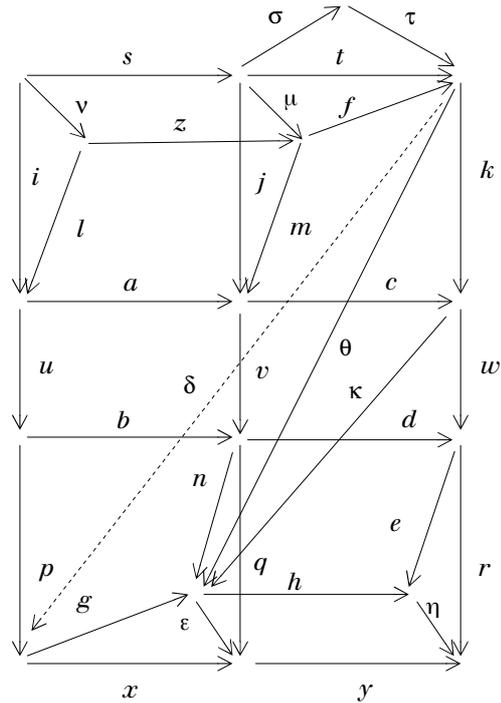

**Figure 9.**

*Construction of diagram 9.* Let $(m,f)$ be the fiber product of $(k,c)$. The square $kf = cm$ is cartesian and since $c$ is an epimorphism, so is $f$ and the square is cocartesian (prop. 5 § 2). Let $z$ be a kernel of $f$. Since $cmz = kfz = 0$ and $a$ is a kernel of $c$, there exists a unique arrow $l$ such that $al = mz$. The square thus built is a kernel of the square built over $m$ and $k$; since $k$ is a monomorphism, this square is cartesian (prop. 5 § 2).

Dualy, one builds the amalgamated sum $(n, g)$ of $(p, b)$. This square is cocartesian and since $b$ is a monomorphism, so is $g$ and the square is cartesian (prop. 2 § 2). Similarly, one builds the cokernel $h$ of $g$ and one completes the square over $h$ and $d$, which is a cokernel of the square built over $g$ and $b$; since $p$ is an epimorphism, this square is cocartesian (prop. 5(c) § 2).

*Proof.* Arrow $nvm$ satisfies $(nvm)z = g(pu)l = 0$ and since $f$ is a cokernel of $z$, there exists a unique arrow $\theta$ such that $nvm = \theta f$. Now $h\theta$ is nul parce que $h\theta f = 0$ and $f$ is an epimorphism. Therefore $\theta$ factorises through the kernel of $h$, that is $g$: there exists a unique arrow $\delta$ such that $\theta = g\delta$. This terminates the construction of $\delta$.

There remains to show that the sequence $(t, \delta)$ is exact or again, since $g$ is a monomorphism, that $(t,\theta)$ is exact; by the duality property $(\delta, x)$ will also be exact. It is already clear that $\theta t$ is nul: $\theta t = nvj = 0$. Let us show that the sequence $(t,\theta)$ is exact.

*Step 1:* Notice that $nva = gpu = 0$ implies that $nv$ factorizes through the cokernel $c$ of $a$: $nv = \kappa c$ for a unique arrow $\kappa$. Moreover, $\kappa k f = \kappa cm = nvm = \theta f$, hence $\kappa k = \theta$ (because $f$ is an epimorphism). Therefore, proving that $(t, \theta)$ is exact reduces to show that $(t, \kappa k)$ exact or again, with the decomposition $t = \tau\sigma$ into an epimorphism followed by a monomorphism, that $\tau$ is a kernel of $\kappa k$.

*Step 2:* The sequence $(va, n)$ is exact. Indeed, in figure 10, the sequence $(u,p)$ is exact, $n(va) = gpu = 0$ and square C is cocartesian by construction. Proposition 1 ensures that the sequence $(va,n)$ is exact. Since $(s, t)$ is exact, $\sigma$ is a cokernel of $s$. Since $c$ is a cokernel of $a$, there exists a unique arrow $\nu$ such that $\nu\sigma = cj$. And $\nu = k\tau$ is a monomorphism, since $k$ and $\tau$ are two monomorphisms.

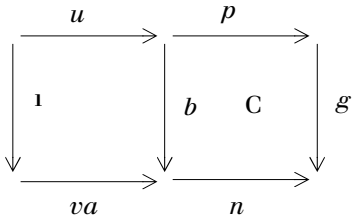

Figure 10.

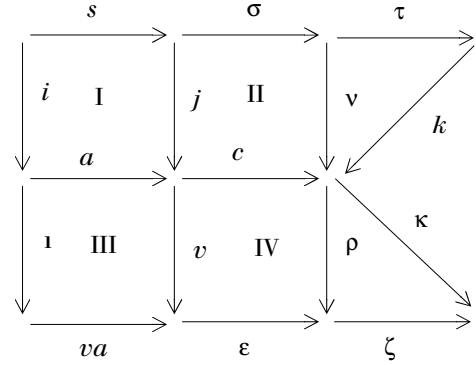

Figure 11.

Decompose $n$ into an epimorphism $\varepsilon$ followed by a monomorphism $\zeta$. Since $(va,n)$ is an exact sequence, $\varepsilon$ is a cokernel of $va$; since $c$ is a cokernel of $a$, there exists a unique arrow $\rho$ such that $\rho c = \varepsilon v$ (see Fig. 11).

*Step 3:* Now, square IV in fig. 11 is a cokernel of square III, and since 1 is an epimorphism, IV is cocartesian (prop. 5(c) § 2). Further, the sequence $(j,v)$ is exact and $\rho\nu = 0$, as can be checked if we precede it with the epimorphism $\sigma$: $\rho\nu\sigma = \varepsilon vj = 0$. From proposition 1, $(\nu,\rho)$ is exact. Since $\nu$ is a monomorphism, it is a kernel of $\rho$, and also of $\zeta\rho = \kappa$ since $\zeta$ is a monomorphism. Proposition 3 of § 1 terminates the proof: $\nu = k\tau$ is a kernel of $\kappa$, therefore $\tau$ is a kernel of $\kappa k = \theta$, QED.